\documentclass[12pt, 14paper,reqno]{amsart}
\usepackage{amsmath,amsfonts,amssymb}
\usepackage[breaklinks]{hyperref}
\usepackage{graphicx}
\usepackage{longtable}
\makeatletter
\@namedef{subjclassname@2020}{%
  \textup{2020} Mathematics Subject Classification}
\makeatother

        \makeatletter
\def\legendre@dash#1#2{\hb@xt@#1{%
  \kern-#2\p@
  \cleaders\hbox{\kern.5\p@
    \vrule\@height.2\p@\@depth.2\p@\@width\p@
    \kern.5\p@}\hfil
  \kern-#2\p@
  }}
\def\@legendre#1#2#3#4#5{\mathopen{}\left(
  \sbox\z@{$\genfrac{}{}{0pt}{#1}{#3#4}{#3#5}$}%
  \dimen@=\wd\z@
  \kern-\p@\vcenter{\box0}\kern-\dimen@\vcenter{\legendre@dash\dimen@{#2}}\kern-\p@
  \right)\mathclose{}}
\newcommand\legendre[2]{\mathchoice
  {\@legendre{0}{1}{}{#1}{#2}}
  {\@legendre{1}{.5}{\vphantom{1}}{#1}{#2}}
  {\@legendre{2}{0}{\vphantom{1}}{#1}{#2}}
  {\@legendre{3}{0}{\vphantom{1}}{#1}{#2}}
}
\def\dlegendre{\@legendre{0}{1}{}}
\def\tlegendre{\@legendre{1}{0.5}{\vphantom{1}}}
\makeatother
\theoremstyle{plain}
\theoremstyle{plain}

\newtheorem{rem}{Remark}[section]
\numberwithin{equation}{section}
\newtheorem{thm}{Theorem}[section]

\newtheorem{cor}{Corollary}[section]
\newtheorem{prop}{Proposition}[section]

\numberwithin{equation}{section}

\usepackage[breaklinks]{hyperref}
\makeatletter
\let\oldHyPsd@CatcodeWarning\HyPsd@CatcodeWarning
\renewcommand{\HyPsd@CatcodeWarning}[1]{
  \ifnum\pdfstrcmp{#1}{math shift}=0    
  \else                                 
    \oldHyPsd@CatcodeWarning{#1}
  \fi
}
\makeatother

\begin{document} 
\title[Diophantine equation]{Solution of certain Diophantine equations in Gaussian integers}
\author{Arkabrata Ghosh}
\address{G. Arkabrata@ SRM University-AP, Mangalgiri Mandal, Neerukonda, Amravati, Andhra Pradesh-522502.}
\email{arka2686@gmail.com}



\maketitle

\section*{Abstract}
In this article, we show that the quartic Diophantine equations $x^4 \pm pqy^4=\pm z^2$ and $ x^4 \pm pq y^4= \pm iz^2$ have only trivial solutions for some primes $p$ and $q$ satisfying conditions $ p \equiv 3 \pmod 8, ~ q \equiv 1 \pmod 8 ~\text{and}~ \displaystyle\legendre{p}{q} = -1$. Here we have found the torsion of the two families of elliptic curves to find the solutions of given Diophantine equations. Moreover, we also calculate the rank of these two families of elliptic curves over the Gaussian field $\mathbb{Q}(i)$.

\section{Introduction}

Integral solution of the Diophantine equation 

\begin{equation}
    \label{eq:1.1}
    ax^4 + by^4= cz^2
\end{equation}

can be easily found in any book on Diophantine equations. By a trivial solution of the equation \eqref{eq:1.1}, we mean to specify the conditions where $x=y=z=0$ or, $a=b=c$ with one of  $x~\text{and}~y$ becomes zero and the square of the other equals $z$. Fermat first studied the Diophantine equation $ x^4 + y^4= z^2$ and using the method of infinite-descent, he proved that no non-trivial integral solution exists to this equation. Later, Hilbert (see \cite{H1894}) extended this result by showing that the equation $x^4 \pm z^2= y^2 $ has only a trivial solution in $\mathbb{Z}[i]$. Now several authors have tried to solve this problem in different aspects. Najman (see \cite{N10}) found all non-trivial solutions of the equation $x^4 \pm y^4 =iz^2$ in $\mathbb{Z}[i]$. Szabo'( see \cite{S04}) solved the Diophantine equations of the form $x^4 + my^4 =z^2$ in $ \mathbb{Z}[i]$ where $ m = \pm 2^n, ~ 0 \neq n \neq 3$. Moreover, in the same article, he proved that if $p \equiv 3 \pmod 8$, then equations of the form $x^4-py^4=z^2$ and $x^2-p^3 y^3 =z^2$ have only trivial solutions in $\mathbb{Z}[i]$ but the equations $x^4+ py^4=z^2$ and $x^2+ p^3 y^3 =z^2$ only have solutions $ (1,1,2)$ and $(2,1,5) $ respectively when $ p =3$. Izadi et al. (see \cite{INB15}) also studied equations of the type \eqref{eq:1.1} over the ring of Gaussian integers. Using the method of elliptic curves, they proved that the families of equations of the form $y^4 \pm p^3 x^4 =z^2$ with $ p \equiv 3 \pmod 8 ~\text{or}~ 3 \pmod {16}$ and $ y^4 \pm px^4= z^2$ with $ 7 ~\text{or}~ 11 \pmod {16}$, have only trivial solutions over $\mathbb{Z}[i]$. Later Ahmadi and Janfada (see \cite{AJ22}) showed that when $ p \equiv 3 \pmod 8$, $q \equiv 1 \pmod 8$ and $\displaystyle\legendre{p}{q} \neq 1$, then the equation  $ x^4 \pm qp^2 y^4= \pm z^2$ and $ x^4 \pm qp^2 y^4 =\pm iz^2$ have only trivial solutions in the ring of Gaussian integers. 

Throughout this article, by a prime element, we mean prime in $\mathbb{Z}$ and we shall refer to primes in $\mathbb{Z}[i]$ as Gaussian primes. Here we take $p$ and $q$ are primes in $\mathbb{Z}$ such that $ p \equiv 3 \pmod 8$, $q \equiv 1 \pmod 8$ and $ \displaystyle\legendre{p}{q} = -1$. Now the main Theorem of this article is as follows.

   \begin{thm}
    \label{Thm:1.1}
    For primes $ p \equiv 3 \pmod 8$ and $ q\equiv 1 \pmod 8$ with $\displaystyle\legendre{p}{q} = -1$, the Diophantine equations 
    \begin{equation}
    \label{eq:1.2}
        x^4  \pm pq y^4= \pm z^2
    \end{equation} 
    
    and 
    \begin{equation}
    \label{eq:1.3}
        x^4 \pm pq y^4  = \pm i z^2
    \end{equation}
 have trivial solutions in $\mathbb{Z}[i]$. 
\end{thm}

The element $ \lambda = 1 +i$ is a Gaussian prime satisfying $ \lambda^4=-4$. The following Corollary gives more equations of the type \eqref{eq:1.1} with trivial solutions over $\mathbb{Z}[i]$.

\begin{cor}
    \label{cor:1.1}
    Let $p,q $ be primes with $ p \equiv 3 \pmod 8$, $ q \equiv 1 \pmod 8$, and $ \displaystyle\legendre{p}{q} = -1$. Then the Diophantine equations 
    \begin{equation}
    \label{eq:1.4}
        x^4 \pm pq y^4= \pm 2^n z^2
    \end{equation}
     and 
     \begin{equation}
     \label{eq:1.5}
         x^4 \pm pq y^4= \pm i2^{n}z^2
     \end{equation}
      have only trivial solutions in $\mathbb{Z}[i]$ for any $ n \in \mathbb{Z}^+$.
\end{cor}

\begin{rem}
    \label{rem:1.2}
    We know $\mathbb{Z}[i]$ is a UFD with the field of fraction $\mathbb{Q}(i)$. By applying Gauss lemma (see \cite{DF04}, Proposition $5$, section $9.3$), we can say that any solution in $\mathbb{Q}(i)$ gives a solution in $\mathbb{Z}[i]$. Hence, throughout this article, we will consider solutions in $ \mathbb{Z}[i] $.
\end{rem}

This article uses elliptic curves to prove the main Theorem \eqref{eq:1.1}. Elliptic curves over the field $\mathbb{Q}(i)$ are not well-known. Here we are interested in elliptic curves of the form $Y^2 =X^3 + kX$ over $ \mathbb{Q}(i)$ where $ k \in \mathbb{Z}$. Selmer-Mordell conjecture showed that the rank of the elliptic curve $ Y^2 = X^3 + pX $, with $p$ prime, has rank $1$ over $\mathbb{Q}$. Bremmer and Cassels(see \cite{BC84}) showed that this conjecture is true for primes $ p \equiv 5 \pmod 8$ and $ p \leq 1000$. 

Consider two families of elliptic curves given by

\begin{equation}
\label{eq:1.6}
    E_{pq}^{+}: ~ Y^2 = X^3 + pqX
\end{equation}

and 

\begin{equation}
    \label{eq:1.7}
    E_{pq}^{-}: ~ Y^2 = X^3 - pqX,
\end{equation}

where $p$ and $q$ are as in Theorem \eqref{Thm:1.1}. At first, we show that the rank of both families of elliptic curves given by equations \eqref{eq:1.6} and \eqref{eq:1.7} over $\mathbb{Q}(i)$ are zero.

\begin{thm}
\label{Thm:1.2}
    For the primes $ p \equiv 3 \pmod 8$ and $ q\equiv 1 \pmod 8$ with $\displaystyle\legendre{p}{q} = -1$, the rank of two families of elliptic curves
given by equations \eqref{eq:1.6} and \eqref{eq:1.7} over $\mathbb{Q}(i) $ are zero.
\end{thm}

Moreover, we also compute the torsion groups of these families of elliptic curves given by \eqref{eq:1.6} and \eqref{eq:1.7} and prove the following Theorem.
\begin{thm}
    \label{Thm:1.3}
     For the primes $ p \equiv 3 \pmod 8$ and $ q\equiv 1 \pmod 8$ with $\displaystyle\legendre{p}{q} = -1$, the torsion group of two families of elliptic curves given by equations\eqref{eq:1.6} and equations \eqref{eq:1.7} over $\mathbb{Q}(i)$ are isomorphic to $\mathbb{Z}/2\mathbb{Z}$.
\end{thm}

\section{Preliminaries}
Let $E$ be an elliptic curve over a field $\mathbb{K}$ of characteristic different from $2$ or $3$ and let $ E(\mathbb{K})$ denote the $\mathbb{K}$-rationals points of $E$ over $\mathbb{K}$. Mordell-Weil Theorem asserts that $  E(\mathbb{K})$ is a finitely generated abelian group and can be represented as

\begin{equation*}
    E(\mathbb{K}) \cong \mathbb{Z}^r \oplus E(\mathbb{K})_{tors}, 
\end{equation*}
where $ r \geq 0$ is called the Mordell-Weil rank of the elliptic curve $E$ over $\mathbb{K}$ and $ E(\mathbb{K})_{tors}$ denotes the torsion sub-group of $ E(\mathbb{K})$ which is a finite abelian group consisting of elements of finite order. When $\mathbb{K}=\mathbb{Q}(i)$, it is well known that there are exactly $16$ possible torsion groups, namely $15$ from Mazur's Theorem and the group $ \mathbb{Z}/4 \mathbb{Z} \times \mathbb{Z}/4 \mathbb{Z}$ ~ (see \cite{N11}). Now to determine the torsion subgroup of the families given by the equation \eqref{eq:1.6},
we need the following Theorem.

\begin{thm}(\textit{Extended Nagell-Lutz Theorem}, see \cite{S09})
    \label{thm:2.1}
    We consider the elliptic curve $Y^2= X^3 + a X + b$ with $a, b \in \mathbb{Z}[i]$ and, let $(X,Y) \in E(\mathbb{Q}(i)) $ be a torsion point. Then, 
    \begin{equation*}
       \begin{cases}
         X,Y \in \mathbb{Z}[i], \\
         \text{Either}~ Y =0, ~\text{or}~ Y^2|(4a^3 + 27b^2). 
        \end{cases}
    \end{equation*}
\end{thm}

To prove Theorem \eqref{Thm:1.2}, we first need to determine the rank of the elliptic curves given by \eqref{eq:1.6} and \eqref{eq:1.7}. So to compute the rank, we need to use the method of $2$-descent. We will describe it briefly here and one can look at \cite{S09} for more details. Suppose that $E : Y^2= X^3 + A X^2 + BX$ is an elliptic curve over $\mathbb{Q}$ and $\overline{E}: Y^2 = X^3 - 2AX^2 + (A^2- 4B) X$ is the corresponding isogenous curve to $E$. Moreover, let $\mathbb{Q}^{*}$ be the multiplicative group of all non-zero rational numbers, and $\mathbb{Q^*}^{2}$ be its subgroup of squares of elements of $ \mathbb{Q^*}$. Now we define the $2$-descent homomorphism $\alpha: E(\mathbb{Q}) \rightarrow \mathbb{Q^*}/\mathbb{Q^*}^{2}$ by

\begin{equation*}
    \alpha(P)= 
               \begin{cases}
                 1 \mod \mathbb{Q^*}^2, ~\text{if}~ P= \mathcal{O}, ~\text{the point of infinity},\\
                 B \mod \mathbb{Q^*}^2, ~\text{if}~ P= (0,0),\\
                 X \mod \mathbb{Q^*}^2, ~\text{if}~ P= (X,Y ) ~\text{with}~ X \neq 0.
                 \end{cases}
\end{equation*}

Similarly, we can define the $2$-decent homomorphism on the isogeneous curve $\overline{E}(\mathbb{Q}) $ curve as follows: $\overline{\alpha}:  \overline{E}(\mathbb{Q}) \rightarrow  \mathbb{Q^*}/\mathbb{Q^*}^{2}$ by 

\begin{equation*}
    \overline{\alpha}(\overline{P})= 
               \begin{cases}
                 1 \mod \mathbb{Q^*}^2, ~\text{if}~ \overline{P}= \mathcal{O}, ~\text{the point of infinity},\\
                 \overline{B} \mod \mathbb{Q^*}^2, ~\text{if}~ \overline{P}= (0,0),\\
                 X \mod \mathbb{Q^*}^2, ~\text{if}~ \overline{P}= (X,Y ) ~\text{with}~ X \neq 0,
                 \end{cases}
\end{equation*}

where $\overline{B}= A^2 - 4B$.

To compute the rank of elliptic curves given by the equation \eqref{eq:1.6} and \eqref{eq:1.7}, we use the following Proposition (see \cite{ST15}, Chapter III, section $6$, page $91$).
\begin{prop}
    \label{prop:2.1} Let $r$ be the rank of $E(\mathbb{Q}) $ and $\alpha$ and $\Bar{\alpha}$ are as above. Then,
    $$
        |\alpha(E(\mathbb{Q}))| |\overline{\alpha}(\overline{E}(\mathbb{Q})|= 2^{r+2}.
    $$
\end{prop}

To compute $  |\alpha(E(\mathbb{Q}))|$ and $ |\overline{\alpha}(\overline{E}(\mathbb{Q})| $, we need to find factorisation of $B$. The following Theorem (see \cite{C07}, Theorem $8.2.9$) precisely says this.

\begin{thm}
    \label{thm:2.2}
    The group $ \alpha(E(\mathbb{Q}))$ is equal to classes modulo squares of $1, B$ and the positive and negative divisors of $A_1$ of $B$ such that the equation
    \begin{equation}
        \label{eq:2.2}
        V^2 = A_1 U^4 +   A U^2 W^2 + (B/A_1) W^4
    \end{equation}
    has a solution $(U, V, W) $ where $U, V ~\text{and}~ W$ are pairwise co-prime  such that $UW \neq 0$ and $gcd(A_1, W)=gcd((B/A_1) , U)=1$. Moreover, there exist $P=(\frac{A_1U^2}{W^2}, \frac{A_1UV}{W^3})$ in $E(\mathbb{Q})$ such that $\alpha(P)=A_1$.
\end{thm}

Now if $\mathbb{K} = \mathbb{Q}(\sqrt{m})$, where $m$ is a square-free integer, we can find the rank of any elliptic curve $E$ over $\mathbb{K}$ by adding ranks of $E$ and its' $m$-twist $E[m]$ over $\mathbb{Q}$. This can be seen from the following result (see \cite{S09}).

\begin{thm}
    \label{thm:2.3}
    Let $\mathbb{K} = \mathbb{Q}(\sqrt{m})$ be a quadratic field, where $m$ is a square-free integer. Let $E: y^2 = x^3 + ax^2 + bx $ be an elliptic curve over $\mathbb{Q}$ and $E[m]: y^2 = x^3 + max^2 + m^2 bx$ be the $m$-twist of $E$. Then

    $$
    rank(E(\mathbb{K}))= rank(E(\mathbb{Q})) + rank(E[m](\mathbb{Q}) ).
    $$
    
\end{thm}

At first, we will prove the Theorem \eqref{Thm:1.2}.

\section{Proof of the Theorem $\texorpdfstring{\eqref{Thm:1.2}}{}$}

\begin{proof}

We prove this result for the family of the elliptic curve given by \eqref{eq:1.7}. The proof for the other family follows a similar logic. According to Proposition \eqref{prop:2.1}, we need to show 

$$
    |\alpha(E_{pq}(\mathbb{Q}))| |\overline{\alpha}(\overline{E_{pq}}(\mathbb{Q}))|= 4.
$$

By the definition of $2$-descent homomorphism of $ \alpha$ and $ \Bar{\alpha}$, it is enough to show

\begin{equation}
\label{eq:3.1}
     |\alpha(E_{pq}(\mathbb{Q}))|= 2 = |\overline{\alpha}(\overline{E_{pq}}(\mathbb{Q}))|.
\end{equation}

We prove the equation \eqref{eq:3.1} in the following steps.

\textit{Step-I} At first, we need to show  that $ |\alpha(E_{pq}(\mathbb{Q}))|= 2$. Here the quartic equation \eqref{eq:2.2} can be written as

\begin{equation}
    \label{eq:3.2}
    V^2 = A_1 U^4 - \frac{pq}{A_1} W^4.
\end{equation}

Therefore, modulo squares, $A_1 \in \{ \pm 1, \pm p, \pm q, \pm pq \}$. From the definition of $\alpha$, we can say $1 $ and $ -pq$ are in the image of $\alpha$. We claim that $Im(\alpha)= \{1, -pq\}$. We discard the other possibilities in the following way.

\textit{Case-I}: If  $ A_1 = -1$, then using the equation \eqref{eq:3.2}, we can say $V^2= -U^4+ pq W^4 $.  Then we have $ (V/U^2)^2  \equiv 1 \pmod p $. By using the Legendre symbol, we can say $\displaystyle\legendre{-1}{p}=1$, that implies $ p \equiv 1 \pmod 4$. This is a contradiction to our hypothesis.

\textit{Case-II}: if $A_1 = p$, then we get $V^2 = p U ^4 - q W^4$. Hence, $ (V/U^2)^2 \equiv p \pmod q$.  It implies   $\displaystyle\legendre{p}{q}=1 $, a contradiction to the hypothesis.

\textit{Case-III}: When $ A_1 = -p$, then  $ V^2 = -pU^4 + q W^4$. Now, $(V/W^2)^2 \equiv q \pmod p$ from which,  we can say $\displaystyle\legendre{q}{p}=1 $. As $ p \equiv  3 \pmod 8$, $q \equiv 1 \pmod 8$ and $\displaystyle\legendre{p}{q}=-1$, the law of quadratic reciprocity gives $ \displaystyle\legendre{q}{p}=-1$. Hence, we arrived at a contradiction.

\textit{Case-IV}: Here we assume $A_1 =q$. If $ q \in Im(\alpha)$, the $ (pq) \cdot (q)=-p \in Im(\alpha)$, which contradicts the Case-III above. The cases $A_1= -q ~\text{and,}~ A_1= -pq$ can be discarded similarly. 

So by combining all different scenarios, it is evident that $Im(\alpha)= \{1, -pq\}$ and, hence $ |\alpha(E_{pq}(\mathbb{Q}))|=2 $.

\textit{Step-II}: We now prove that $ |\overline{\alpha}(\overline{E_{pq}}(\mathbb{Q}))|= 2$. The quartic equation \eqref{eq:2.2} can be written as

\begin{equation}
    \label{eq:3.3}
    \overline{V}^2 =  \overline{A_1}~ \overline{U}^4 + \dfrac{4pq}{\overline{A_1}} \overline{W}^4.
\end{equation}

So, modulo squares,

$$
\overline{A_1} \in \{ \pm 1, \pm 2,  \pm p, \pm q, \pm pq, \pm 2p, \pm 2q, \pm 2pq \}. 
$$

By the definition of $ \overline{\alpha}$, $1 ~\text{and}~ pq $ are in the image of $\overline{\alpha}$. We will show no other elements in the image of $\overline{\alpha}$ in the following way: We consider all other elements of $\overline{A}$ and arrive at a contradiction in each case. 

\textit{Case-I}: $\overline{A_1}=-1$; Using equation \eqref{eq:3.3}, we get $  \overline{V}^2 = -(\overline{U}^4 + 4pq\overline{W}^4)$, a clear contradiction. Proceeding similarly, we can say $\overline{A}$ can't be $ -p, -q , -2, -pq, -2p, -2q, ~\text{and}~ -2pq$.

\textit{Case-II}:  $\overline{A_1}=p$; Hence, we have $  \overline{V}^2 = (p\overline{U}^4 + 4q\overline{W}^4)$. It implies $  \bigg(\overline{V}/\overline{W}^2\bigg)^2 \equiv q \pmod p$, which says  $\displaystyle\legendre{p}{q}=1 $, a contradiction to the hypothesis.

\textit{Case-III}: $\overline{A_1}=q$; $  \overline{V}^2 = (q\overline{U}^4 + 4p\overline{W}^4)$. Then, $ \bigg(\overline{V}/\overline{U}^2\bigg)^2 \equiv q \pmod p$. As $ p \equiv  3 \pmod 8$, $q \equiv 1 \pmod 8$ and $\displaystyle\legendre{p}{q}=-1$, the law of quadratic reciprocity gives $ \displaystyle\legendre{q}{p}=-1$. So we get a contradiction.

\textit{Case-IV}: $\overline{A_1}=2pq$; $  \overline{V}^2 = (2pq\overline{U}^4 + 2\overline{W}^4)$. Now, $ \bigg(\overline{V}/\overline{W}^2\bigg)^2 \equiv 2 \pmod p$ from which, we get $\displaystyle\legendre{2}{p}=1$. This gives $ p \equiv 1 ~\text{or}~ 7 \pmod 8$, which is impossible.

\textit{Case-V}: $\overline{A_1}=2$; Now if $ 2 \in Im(\overline{{\alpha}})$, then $ (pq) \cdot 2= 2pq \in Im(\overline{{\alpha}}) $, which is a contradiction according to case-IV from above. Using a similar technique, we can show that $ 2p ~\text{and}~ 2q$ can't be values of $\overline{A_1}$.

Now combining all cases discussed in Step-II, we get $  |\overline{\alpha}(\overline{E_{pq}}(\mathbb{Q}))|= 2$. Hence, the equation \eqref{eq:3.1} is satisfied and we get that rank$(E(\mathbb{Q}))=0$.

If we take $ m=-1$, then $ m$-twist of $E_{pq}$ is itself, and hence by applying the Theorem \eqref{thm:2.3}, we can conclude the proof of Theorem \eqref{Thm:1.2}.

\end{proof}


\section{Proof of the Theorem $\texorpdfstring{\eqref{Thm:1.3}}{}$}

\begin{proof}
Let us assume $(X, Y) $ is a torsion point of the elliptic curve $E_{pq}^{+}$ given by the equation \eqref{eq:1.6}. From the Theorem \eqref{thm:2.1}, we know that $X, Y \in \mathbb{Z}[i]$, and moreover, either $Y =0$ or $ Y^2 \mid 4p^3 q^3$. Now, $Y=0$ gives the $2$-torsion point $(0,0)$. We claim that $Y^2 \mid 4p^3 q^3$ gives no other torsion points. Suppose our assumption is wrong. Then $Y^2 \mid 4p^3 q^3$ implies $Y^2= lq^u p^v$ where 
\begin{equation}
\label{eq:3.4}
    l \in \{\pm 1, \pm 2i, \pm 4 \}, ~ u,v \in \{0,2\}.
\end{equation}
At first, we will ignore the sign of $l$  and prove the contradiction in all twelve cases. The opposite signs of $l$ can be discarded similarly. Throughout this proof, we take $ \lambda= (1+i)$.
  
   \textit{Scenario-I}: Let $Y^2 =2i =\lambda^2 = X^3 + pqX$. The only Gaussian prime factor of $X$ is $ \lambda$ and we can take $X=\lambda^{k}$ for some $ k \geq 1$. Then, we have 
    \begin{equation}
    \label{eq:3.5} 
        \lambda^2 = \lambda^{3k} + pq \lambda^{k}.
    \end{equation}
     . Now we will consider two cases.

    \textit{Case-I}: At first, we take $k=1$. Putting this value in the equation \eqref{eq:3.5}, we get $\lambda^2 = \lambda^3 + pq \lambda$ which after simplification, gives $ \lambda- \lambda^2 = pq$. As $ \lambda- \lambda^2 = 1 -i  $, it can't be equal to $pq $. So in this case, equation \eqref{eq:3.5} does not hold.

\textit{Case-II}: Now we consider scenario when $ k \geq 2$. This scenario can be again divided into two sub-cases. 

    \textit{Sub-Case (a)}: Now we will deal with the situation when $k=2$. Then equation \eqref{eq:3.5} can be written as $ \lambda^2 = \lambda^2(\lambda + pq)$, which after simplification, yields $ pq= -i$, a contradiction.

    \textit{Sub-Case~ (b)}: Finally, we take into account the scenario when $k \geq 3$. From the equation \eqref{eq:3.5}, we have $ \lambda^2 = \lambda^{3k} + pq \lambda^{k}$, and after simplifying it, we get $ 1= \lambda(\lambda^{3k-3} + pq \lambda^{k-3} )$. But from this expression, we can say $ \lambda$ is a unit in $ \mathbb{Z}[i]$. It implies $ \lambda \in \mathbb{Z}[i]^*= \{ \pm 1, \pm i\}$. As $ \lambda = 1+i$, we arrived at a contradiction.

    Hence. the equation \eqref{eq:3.5} does not hold and we arrived at a contradiction. Similarly, we can also show that $ Y^2 = 4 =\lambda^4$ is impossible.

   \textit{Scenario-II}: Let $Y^2 = q^2 = X^3 + pq X$. Then $q \mid X$. Moreover, let us assume that $m$ is the highest index of $q $ dividing $X$. It implies $X=q^{m} X_0$ for some $ m\geq 1$ and $ q $ does not divide $X_0$. Then $q^2= q^{3m}X_0^3 + pq^{m+1}X_0$ and by cancelling $q^2$, we get the equation
    \begin{equation}
        \label{eq:3.6}
        1= q^{3m-2} X_0^3 + pq^{m-1}X_0.
    \end{equation}
    We claim that the equation \eqref{eq:3.5} has no solution in $\mathbb{Z}[i]$. \\
   
    \textbf{Proof of the Claim}:
    \begin{proof}
        
    We will consider two different cases  which are as follows:

    \textit{Case-I}: At first, we assume $ m \geq 2$. Then  $  3m -2 > 1$ and $ m-1 \geq 1$. Hence, considering the \eqref{eq:3.5} modulo $q$, we get $ 1 \equiv 0 \pmod q$, which is a contradiction.

   \textit{Case-II}: Now we consider $m=1$. Then the equation \eqref{eq:3.5} can be written as

   \begin{equation*}
       1= X_0 ( p+ qX_0^2).
   \end{equation*}
    It implies that both $X_0$ and $ ( p + qX_0^2)$ are units in $\mathbb{Z}[i]$. As $ \mathbb{Z}[i]^*=  \{ \pm 1, \pm i \}$, both $ X_0$ and $ ( p + qX_0^2)$ belongs to $ \mathbb{Z}[i]^*$. If $X_0 = \pm 1$, then $ p + qX_0^2 = p+q$ which clearly not belongs to $ \mathbb{Z}[i]^*$. Now when $ X_0 = \pm i$, then $ p + qX_0^2 = p-q$ and clearly $ (p-q) \neq  \pm i$. As both $p$ and $q$ are odd primes, $ p-q = \pm 1$ is impossible. \\
    This concludes the proof of our claim and we can say that equation \eqref{eq:3.5} has no solutions in $\mathbb{Z}[i]$.
    \end{proof}

    \textit{Scenario-III}: For $ Y^2= 2q^2 i=\lambda^2 q^2 = X^3 + pqX$. we again have that $ q|X $.  Here again, we take $m$ to be the highest index of  $q$ dividing $X$ and we can assume $X=q^{m} X_0$ for some $ m\geq 1, ~ q \nmid X_0$. Then we have $ \lambda^2 q^2 = q^{3m} X_0^3 + q^{m+1} X_0$ which after simplification gives 
    $\lambda^2= q^{3m-2} X_0 ^3 + p q^{m-1} X_0$. Now by arguing similarly to \textit{Scenario-II}, we will get a contradiction. The case $ Y^2= 4q^2= \lambda^{4}q^2$ can be discarded similarly.

    \textit{Scenario-IV}: Finally, let $Y^2= lq^{u} p^{v}= X^3 + pqX$ where $l, u, ~\text{and}~ v$ are as described as in equation \eqref{eq:3.4}. From the equation, it is evident that $ p \mid X$.
    Hence, 
    \begin{equation}
    \label{eq:3.7}
        l q^{u} p^{v} = p^{3m'}X_0^3 + p^{m'+1}q X_0,
    \end{equation}
    where $ m' \geq 1$ is the highest index of $p $ in $X$, where $ p \nmid X_0$. Now we will consider two different scenarios. 
    
    \textit{Case-I}: We  assume $ v= min \{v, m'+1 \}$. Hence, by cancelling $p^{v}$ from equation \eqref{eq:3.7}, we get 

    \begin{equation*}
        l q^{u}= p^{3m'-v} X_0^3 + p^{m'+1-u} q X_0.
    \end{equation*}

    From the above equation, it is clear that $ p \mid lq^{u}$, is impossible as $ p $ divides neither $l$ nor $q$.

\textit{Case-II}: Now we take $ m'+1 = min \{v, m'+1 \}$.  Now, similarly by  cancelling $p^{m'+1}$ from equation \eqref{eq:3.7}, we get 

    \begin{equation*}
        l q^{u} p^{v -(m'+1)} = p^{2m'-1} X_0^3 +  q X_0.
    \end{equation*}

    It implies that $ p \mid qX_0$,  a clear contradiction as $ p \nmid q$ and $ p \nmid X_0$.

    So now by combining all the scenarios discussed above, we can say that our assumption is wrong and $ Y^2 \nmid 4 q^3 p^3$. Hence the torsion sub-group of the elliptic curve defined by the equation \eqref{eq:1.6} is isomorphic to $\mathbb{Z}/2\mathbb{Z}$. By proceeding similarly as above, we can also prove that the torsion sub-group of the elliptic curve given by equation \eqref{eq:1.7} is also isomorphic to $\mathbb{Z}/2\mathbb{Z}$. This concludes the proof of Theorem \eqref{Thm:1.2}.
\end{proof}

Finally, we will prove the main Theorem \eqref{Thm:1.1}

\section{Proof of the Theorem $\texorpdfstring{\eqref{Thm:1.1}}{}$}

\begin{proof}
    We prove the results for the equations \eqref{eq:1.2} and \eqref{eq:1.3} with only positive signs on the right-hand side as the negative side has similar proofs.\\
    Let us assume that $(x,y,z)$ is a non-trivial solution of the equation \eqref{eq:1.2}. Dividing both sides by $y^4$ and taking $x/y ~ = ~ a$ and $z/y^2 ~ = ~ b$, we get $ a^4  \pm pq = b^2$. By taking $X=a^2$, we have two equations

    \begin{equation*}
        X = a^2, ~ X^2 \pm pq = b^2.
    \end{equation*}

Now, multiplying these two equations and assuming $ Y ~= ~ ab$, we obtain a torsion point $(X,Y) ~ = (ab, a^2) \neq (0,0) $ on the elliptic curve $ Y^2 = X^3 \pm pq X$, clear contradiction to Theorem \eqref{Thm:1.3}.

Now, we will consider the equation \eqref{eq:1.3}. We may assume $(x,y,z) $ as a non-trivial solution to the equation \eqref{eq:1.3} and divide both sides of it by $y^4$. Now similarly to the above scenario, we take  $ x/y ~ = ~ a$ and $z/y^2 ~ = ~ b$.  Then we get $ a^4  \pm pq = ib^2$. Now by taking $ -iX=a^2$, we obtain following equations

\begin{equation*}
        -iX = a^2, ~ -X^2 \pm pq = ib^2.
    \end{equation*}
Proceeding similarly as in the case of equation \eqref{eq:1.3}, we get a torsion point $ \neq (0,0)$ on the elliptic curve $ Y^2 =  X^3 \pm pq X$, leading again to a contradiction.

\end{proof}

Now we will end this article with the proof of  Corollary \eqref{cor:1.1}

\section{Proof of the Corollary $\texorpdfstring{\eqref{cor:1.1}}{}$}

\begin{proof}

If we consider $n=2k $ and $n=2k+1$ in the equation \eqref{eq:1.4}, we get , respectively, the equations

\begin{equation*}
    x^4 \pm pqy^4= \pm (2^{k} z)^2, ~ ~ x^4 \pm pqy^4 = \pm i(i\lambda^2 2^{k} z)^2.
\end{equation*}

These two equations are clearly in the form of the equations \eqref{eq:1.2} and \eqref{eq:1.3} of Theorem \eqref{Thm:1.1}. 

Similarly, if we consider the equation \eqref{eq:1.5} and take $n=2k$ and $n= 2k+1$ respectively, then we get

\begin{equation*}
    x^4 \pm pqy^4= \pm i(2^{k} z)^2, ~ ~ x^4 \pm pqy^4 = \pm(\lambda^2 2^k z)^2.
\end{equation*}

Hence, we get the desired conclusion from the Theorem \eqref{Thm:1.1}.

\end{proof}

\section*{Funding and Conflict of Interests/Competing Interests} The author has no financial or non-financial interests to disclose that are directly or indirectly related to the work. The author has no funding sources to report.

\section*{Data availability statement} No outside data was used to prepare this manuscript.

\section*{Acknowledgement} The author is grateful to the anonymous referee for his/her valuable suggestions for the improvement of this manuscript.

\end{document}